\documentclass[reqno,centertags, 12pt, draft]{article}

\usepackage{amsmath,amsthm,amscd,amssymb}

\usepackage{latexsym}

\usepackage{mathrsfs}

\newcommand{\bbN}{{\mathbb{N}}}

\newcommand{\bbR}{{\mathbb{R}}}

\newcommand{\bbZ}{{\mathbb{Z}}}

\newcommand{\bbC}{{\mathbb{C}}}

\newcommand{\calL}{{\mathcal L}}

\newcommand{\scrV}{{\mathscr V}}

\newcommand{\ip}[2]{\left( #1 \left | #2 \right. \right)}

\newcommand{\no}{\nonumber}

\newcommand{\ti}{\tilde  }

\newcommand{\beq}{\begin{equation}}

\newcommand{\eeq}{\end{equation}}

\newcommand{\ba}{\begin{align}}

\newcommand{\ea}{\end{align}}

\DeclareMathOperator{\Ima}{Im}

\allowdisplaybreaks

\numberwithin{equation}{section}

\newtheorem{theorem}{Theorem}[section]

\newtheorem{proposition}[theorem]{Proposition}

\theoremstyle{definition}

\theoremstyle{remark}

\newtheorem*{remarks}{Remarks}

\begin{document}
\title{Localization for the Anderson Model on Trees with Finite 
Dimensions}
\author{Jonathan Breuer \\
\footnotesize Institute of Mathematics, \\
\footnotesize The Hebrew University of Jerusalem, \\
\footnotesize 91904 Jerusalem, \\
\footnotesize Israel. \\
\footnotesize Email: jbreuer@math.huji.ac.il}
\maketitle
\begin{abstract}
We introduce a family of trees that interpolate between the Bethe lattice and $\bbZ$. We prove complete localization for 
the Anderson model on any member of that family.
\end{abstract}

\section{Introduction}
\sloppy
The purpose of this paper is to study the spectral properties of the Anderson model on a family of graphs which interpolate 
in a certain sense between the Bethe lattice and $\bbZ$. The Bethe lattice can be regarded as an infinite dimensional 
graph because of the exponential growth (in $r$) of the volume of the ball of radius $r$ around the root (which is 
connected to the fact that a significant part of this volume is concentrated on the boundary of that ball). $\bbZ$ is, 
of-course, a one-dimensional graph. For both these domains, the Anderson model has been extensively studied. 
For the one-dimensional case it is known that the spectrum is pure point 
(with exponentially decaying eigenfunctions) for all energies and any degree of disorder 
(see, e.g., \cite{carmona} and references therein). On the other hand, 
absolutely continuous spectrum is known to occur, on the Bethe lattice, in the weak disorder regime 
\cite{aiz1, froese-bethe, klein-bethe}. 
We present here a family of trees whose members are all finite dimensional in a natural sense and which have the Bethe 
lattice and $\bbZ$ as extreme cases. In this setting, we shall prove localization 
for the Anderson model for all energies and any degree of disorder whenever the dimension is finite.

In order to describe the objects at the focus of our attention, we need some terminology. 
By \emph{a rooted tree}, $\Gamma$, we mean a tree graph that has a special vertex designated by the letter $O$. We use 
$\mathscr{V}(\Gamma)$ to 
denote the set of vertices of $\Gamma$. For any two vertices, $x, y \in \mathscr{V}(\Gamma)$ it is
possible to define the distance between $x$ and $y$, $d(x,y)=d(y,x)$, as the number
of edges of the unique path of minimal length connecting them. 
These notions allow us to define a natural direction on the tree: For any vertex on a rooted tree, the backward direction
is the direction pointing towards the root. Any other direction we call \emph{forward}. More precisely, for  $x \in 
\mathscr{V}(\Gamma)$, we say that $y$ is a \emph{forward} neighbor of 
$x$ if $d(x,y)=1$ and $d(y,O)>d(x,O)$. In this case, we shall say that $x$ is a \emph{backward} neighbor of $y$. 

The trees we construct are parametrized by a natural number $k \geq 2$ and a real number $\gamma \geq 1$. Roughly speaking, 
they are obtained by taking the Bethe lattice of coordination number $k$, and extending its edges at an exponential 
rate (determined by $\gamma$). This is done by replacing the edges at a distance $n$ from the root by a segment of $\bbZ$ 
of length $[\gamma^n]$ (where $[\ \cdot \ ]$ for a real number denotes its integer part).

More precisely, let $k \geq 2$ be a natural number and $\gamma\geq 1$ be a 
real number. We define the rooted tree $\Gamma_{k,\gamma}$ as follows: 
Let $\mathscr{S}_{k,\gamma} \subseteq \mathscr{V}(\Gamma_{k,\gamma})$ be the set of vertices of $\Gamma_{k,\gamma}$ whose elements are 
the root $O$, and all vertices at a 
distance $\sum_{j=1}^N [\gamma^j]$ from $O$ (for any $N\in \bbN$). Now, $\Gamma_{k,\gamma}$ is defined by the fact that vertices 
belonging to 
$\mathscr{S}_{k,\gamma}$ have $k$ forward neighbors. All other vertices have one forward nearest neighbor (see Figure 1). 
We call the elements of 
$\mathscr{S}_{k,\gamma}$ \emph{junctions}. 
It is easy to see that by taking $\gamma=1$ we get the Bethe lattice of coordination number $k$. 
On the other hand, $\bbZ$ can be viewed as corresponding to the case  $k=2$, 
$\gamma=\infty$. In this sense, the family $\{\Gamma_{k,\gamma}\}_{k \geq 2, \gamma\geq 1}$ 
interpolates between the Bethe lattice and $\bbZ$.

\setlength{\unitlength}{0.7cm}
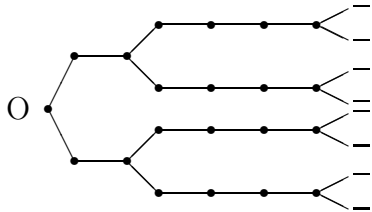
\begin{figure}[!hbp]
\begin{picture}(5,4)(0,-2)
\put(2.2,-0.2) O
\put(3,0){\circle*{0.15}}
\put(3,0){\line(1,2){0.5}}
\put(3.5,1){\circle*{0.15}}
\put(3,0){\line(1,-2){0.5}}
\put(3.5,-1){\circle*{0.15}}
\put(3.5,1){\line(1,0){1}}
\put(4.5,1){\circle*{0.15}}
\put(3.5,-1){\line(1,0){1}}
\put(4.5,-1){\circle*{0.15}}
\put(4.5,-1){\line(1,1){0.6}}
\put(4.5,-1){\line(1,-1){0.6}}
\put(4.5,1){\line(1,1){0.6}}
\put(4.5,1){\line(1,-1){0.6}}
\put(5.1,1.6){\line(1,0){1}}
\put(6.1,1.6){\circle*{0.15}}
\put(6.1,1.6){\line(1,0){1}}
\put(7.1,1.6){\circle*{0.15}}
\put(7.1,1.6){\line(1,0){1}}
\put(8.1,1.6){\circle*{0.15}}
\put(5.1,0.4){\line(1,0){1}}
\put(6.1,0.4){\circle*{0.15}}
\put(6.1,0.4){\line(1,0){1}}
\put(7.1,0.4){\circle*{0.15}}
\put(7.1,0.4){\line(1,0){1}}
\put(8.1,0.4){\circle*{0.15}}
\put(5.1,-0.4){\line(1,0){1}}
\put(6.1,-0.4){\circle*{0.15}}
\put(6.1,-0.4){\line(1,0){1}}
\put(7.1,-0.4){\circle*{0.15}}
\put(7.1,-0.4){\line(1,0){1}}
\put(8.1,-0.4){\circle*{0.15}}
\put(5.1,-1.6){\line(1,0){1}}
\put(6.1,-1.6){\circle*{0.15}}
\put(6.1,-1.6){\line(1,0){1}}
\put(7.1,-1.6){\circle*{0.15}}
\put(7.1,-1.6){\line(1,0){1}}
\put(8.1,-1.6){\circle*{0.15}}
\put(5.1,1.6){\circle*{0.15}}
\put(5.1,0.4){\circle*{0.15}}
\put(5.1,-1.6){\circle*{0.15}}
\put(5.1,-0.4){\circle*{0.15}}
\put(8.1,-1.6){\line(2,1){0.6}}
\put(8.1,-1.6){\line(2,-1){0.6}}
\put(8.1,-0.4){\line(2,1){0.6}}
\put(8.1,-0.4){\line(2,-1){0.6}}
\put(8.1,1.6){\line(2,1){0.6}}
\put(8.1,1.6){\line(2,-1){0.6}}
\put(8.1,0.4){\line(2,1){0.6}}
\put(8.1,0.4){\line(2,-1){0.6}}
\put(8.8,-1.9){\line(1,0){0.3}}
\put(8.8,-1.25){\line(1,0){0.3}}
\put(8.8,-0.7){\line(1,0){0.3}}
\put(8.8,-0.05){\line(1,0){0.3}}
\put(8.8,1.95){\line(1,0){0.3}}
\put(8.8,1.3){\line(1,0){0.3}}
\put(8.8,0.75){\line(1,0){0.3}}
\put(8.8,0.15){\line(1,0){0.3}}
\end{picture}
\caption{A neighborhood of the root for $\Gamma_{2,2}$.}
\end{figure}

A straightforward computation shows:
\begin{proposition} \label{dimension}
Fix $\gamma>1$ and $\bbN \ni k \geq 2$. Let $\Gamma=\Gamma_{k,\gamma}$ and let 
$B_\Gamma(r) = \{x \in \mathscr{V}(\Gamma) \mid d(x,O) \leq r \}$. Then
\beq \label{dimensions}
\lim_{r \rightarrow \infty} \frac{\log \#B_\Gamma(r)}{\log r}=\frac{\log \gamma k}{\log \gamma}
=1+\frac{\log k}{\log \gamma}.
\eeq
\end{proposition}
Below, we shall refer to the quantity $\frac{\log \gamma k}{\log \gamma}$ as the \emph{dimension} of $\Gamma_{k,\gamma}$.

Since we are dealing with non-regular trees (namely, the number of nearest neighbors is not constant), there are two 
choices for the Laplacian:
\beq \label{lap1}
(\ti{\Delta} f)(x)=\sum_{y : d(x,y)=1}f(y),
\eeq
and
\beq \label{lap2}
(\Delta f)(x)=\sum_{y : d(x,y)=1 }f(y)-\#\{y : d(x,y)=1 \}\cdot f(x)
\eeq
where $\#A$, for a finite set $A$, is the number of elements in $A$. 
Both operators are bounded and self-adjoint on $\Gamma_{k,\gamma}$ for any $k$ and $\gamma$. Theorem \ref{localization} below
holds, as stated, both for $\Delta$ and $\ti{\Delta}$. Moreover, we shall give a proof that goes through for both cases. 
To avoid encumbrance, we shall use the notation $\Delta$ with the understanding that all statements hold for \eqref{lap1} 
as well as for \eqref{lap2}. 

Let $\Gamma$ be a tree and let $\{V_\omega(x)\}_{x \in \mathscr{V}(\Gamma)}$ be a family of i.i.d.\ random variables with 
common probability distribution $d\rho$. For 
any $\omega$, let $V_\omega$ stand for the corresponding multiplication operator defined over $\ell^2(\Gamma) \equiv 
\ell^2 \left( \mathscr{V}(\Gamma) \right)$ by
\beq \no
(V_\omega f)(x)=V_\omega(x)f(x). 
\eeq
For $\lambda>0$ we refer to the 
family of operators
\beq \no
H_{\omega,\lambda}=\Delta+\lambda V_\omega
\eeq
as the Anderson model with coupling constant $\lambda$. For $\Gamma=\bbZ$, this model is known to exhibit almost sure pure point 
spectrum with exponentially decaying eigenfunctions, for any range of energies and any value of the coupling constant, for any 
probability distribution $d\rho$, either having an absolutely continuous 
component, or having some finite 
moment \cite{carmona}. For the case of the Bethe lattice, on the other hand, it is known 
\cite{aiz1, froese-bethe, klein-bethe} that the Anderson model exhibits absolutely 
continuous spectrum for small values of $\lambda$ (for $d\rho$ satisfying certain regularity 
conditions). 

We shall assume throughout that 

(i) $d\rho$ has a bounded density, namely
\beq \label{Minami-condition-1}
d\rho(\xi)=\ti{\rho}(\xi)d\xi
\eeq
with
\beq \label{Minami-condition-2}
||\ti{\rho}||_\infty<\infty.
\eeq

(ii) \beq \label{Minami-condition-3}
\int |\xi|^\eta \ti{\rho}(\xi) d\xi<\infty \textrm{ for some } \eta>0.
\eeq

Our main result is
\begin{theorem} \label{localization}
Let 
\beq \no
H_{\omega,\lambda}=\Delta+\lambda V_\omega
\eeq
be the Anderson model on $\Gamma=\Gamma_{k,\gamma}$ for some $k \geq 2$ and $\gamma>1$. Assume that $d\rho$ satisfies requirements 
(i)-(ii) above. Then, for any $\lambda>0$ and almost every realization of 
$V_\omega$, $H_{\omega,\lambda}$ has only pure point spectrum and the 
corresponding eigenfunctions decay exponentially.
\end{theorem}
\begin{remarks}
1. For a function $f$ defined on $\mathscr{V}(\Gamma)$, we say that $f$ decays exponentially if there exist positive constants $A,C$ 
such 
that,
\beq \no
f(x) \leq A e^{-C|x|}
\eeq
where $|x|=d(x,O)$.

2. We note that the technical requirement \eqref{Minami-condition-3} is also present in the 
proof of localization for the Anderson-Bernoulli model in one-dimension \cite{ckm}. We, 
however, assume in addition the absolute continuity of $d\rho$ \eqref{Minami-condition-1}, so the question of 
localization for the Anderson-Bernoulli model on $\Gamma_{k,\gamma}$ is still open.
\end{remarks}

The proof of Theorem \ref{localization} relies on the fact that as long as 
$\gamma>1$, $\Gamma_{k,\gamma}$ contains arbitrarily long 
one-dimensional segments. We call trees with this property \emph{sparse}. Applying ideas of the finite-volume method developed by 
Aizenman, Schenker, Friedrich and Hundertmark in \cite{aiz}, we use a priori bounds that are known for the one-dimensional case in 
order to get exponential decay of fractional moments of the Green function. Since $\Gamma_{k,\gamma}$ has finite dimensions in 
the sense of \eqref{dimensions}, this implies localization.

We note that the behavior of the Anderson model on these sparse, finite dimensional trees is drastically different from the expected 
behavior on $\bbZ^d$, where some absolutely continuous spectrum is believed to exist in the weak coupling regime. Such a difference is 
also manifest in the spectral properties of the Laplacian. The papers \cite{breuer-CMP, breuer-mol} are devoted to the 
spectral analysis of $\Delta$ on sparse trees. Examples are constructed, in these papers, of sparse trees where $\Delta$ has singular 
spectral measures. In particular, it is shown in \cite{breuer-mol} that generically, in a certain probabilistic sense, the 
finite dimensional trees discussed here have singular spectrum and some even exhibit a component of dense point type.

We are grateful to Michael Aizenman, Nir Avni, Vojkan Jak\u si\'c, Yoram Last, Barry Simon and Simone Warzel for useful discussions. 
We also wish to thank 
Michael Aizenman for the hospitality of Princeton where this work was done.

This research was supported in part by THE ISRAEL SCIENCE FOUNDATION (grant no. 188/02) and by Grant no. \mbox{2002068} from the 
United States-Israel Binational Science Foundation (BSF), Jerusalem, Israel.
 
\section{Proof of Theorem \ref{localization}}

Fix $\bbN \ni k \geq 2$, $\gamma>1$ and $\lambda>0$. To streamline the notation, let $\Gamma=\Gamma_{k,\gamma}$ and 
$H_\omega=H_{\omega,\lambda}$. We also use the shorthand $|x| \equiv d(x,O)$.
As mentioned earlier, for any $x,y \in \scrV(\Gamma)$, there is a unique path of minimal length connecting them to each 
other. This is a finite subgraph of $\Gamma$ which can also be embedded in $\bbZ$. 
We denote this graph by $\mathcal{L}(x,y)$.

For $x \in \mathscr{V}(\Gamma)$, the spectral measure $\mu_x$ is defined by the equation
\beq \no
\int_\bbR \frac{d\mu_x}{x-z}=\ip{\delta_x}{(H_{\lambda,\omega}-z)^{-1} \delta_x} \qquad z \in \bbC \setminus \bbR
\eeq
where $\delta_x$ is the delta function at $x$ and $\ip{f}{g}$ stands for the inner product in $\ell^2(\Gamma)$.
We shall prove Theorem \ref{localization} by showing that, with probability one, $\mu_x$ is pure point for any 
$x \in \mathscr{V}(\Gamma)$. Since $\{\delta_x\}_{x \in \mathscr{V}(\Gamma)}$ is an orthogonal 
basis for $\ell^2(\Gamma)$, the theorem is immediately implied.

We want to apply ideas of Aizenman et al.\ \cite{aiz} to our setting. 
In particular, we will study $H_\omega$ restricted to 
finite regions of $\Gamma$. For any such region, $\Omega$, we denote 
by $\Theta(\Omega)$ the set of nearest-neighbor 
bonds reaching out of $\Omega$, that is,
\beq \label{Theta(Omega)}
\Theta(\Omega)=\{ (x,x')\in \scrV(\Gamma) \times \scrV(\Gamma) \mid x \in \scrV(\Omega),\ x' \in 
\scrV (\Gamma) \setminus \scrV(\Omega),\ d(x,x')=1 \}.
\eeq
We further denote by $\Omega^+$ the region containing the vertices within distance $1$ from $\Omega$, and by 
$\mathscr{B}(\Omega)$ the boundary of $\Omega$, that is,
\beq \label{B(Omega)}
\mathscr{B}(\Omega)=\{x \in \scrV(\Omega) \mid \exists x' \in \scrV(\Gamma) \setminus \scrV(\Omega) 
\textrm{ s.t. } d(x,x')=1 \}.
\eeq 

We let $\Delta_{\Omega}$ be the operator obtained by
deleting the hopping terms corresponding to $\Theta(\Omega)$ 
(following \cite{aiz} we shall call the off-diagonal matrix elements 
of $\Delta$ \emph{hopping terms}), so that the 
restriction of $\Delta_{\Omega}$ to 
$\ell^2(\Omega)$ is just the finite volume Laplacian with Dirichlet 
boundary conditions on the boundary of $\Omega$. With this we may define the restriction of $H_\omega$ as:
\beq \no
H_{\Omega,\omega}=\Delta_\Omega+V_\omega.
\eeq
For $H_{\Omega,\omega}$ as well, the restriction to $\ell^2(\Omega)$ equals the finite volume operator with Dirichlet boundary
conditions on the boundary of $\Omega$. 

Another kind of restriction we will consider is $H^{\Omega}_\omega$ ($\Delta^\Omega$) which 
is the operator one gets from $H_\omega$ ($\Delta$) by deleting
all hopping terms outside of $\Omega$. The restriction of this 
operator to $\ell^2(\Omega)$ is again just the finite volume
operator with Dirichlet boundary conditions on $\Omega$.

We want to obtain decay of fractional moments of the Green function 
for the operators mentioned above. 
This is the function:
\beq \label{Green}
G_\omega(x,y;z) \equiv \ip{\delta_{x}} {G_\omega(z) \delta_{y}} \equiv
\ip{\delta_{x}}{(H_\omega-z)^{-1} \delta_{y}}
\eeq
defined for any $z$ in the resolvent set of $H_\omega$ and in particular 
for any $z \in \bbC \setminus \bbR$. We also use 
\beq \label{Green_restricted1}
G_{\Omega;\omega}(x,y;z) \equiv \ip{\delta_{x}}{G_{\Omega;\omega}(z) \delta_{y}} \equiv
\ip{\delta_{x}}{(H_{\Omega;\omega}-z)^{-1} \delta_{y}},
\eeq 
and
\beq \label{Green_restricted2}
G^\Omega_\omega(x,y;z) \equiv \ip{\delta_{x}}{G^\Omega_\omega(z) \delta_{y}} \equiv
\ip{\delta_{x}}{(H^\Omega_\omega-z)^{-1} \delta_{y}}.
\eeq

We note that $G_{\Omega;\omega}$ is a direct sum of operators, one 
corresponding to $\Omega$ and the other corresponding
to $\Gamma \setminus \Omega$ so that if $x \in \scrV(\Omega)$ and $y \in
\scrV(\Gamma) \setminus \scrV(\Omega)$, then
\beq \no
\ip{\delta_x }{G_{\Omega,\omega}(z) \delta_y} =0. 
\eeq
The same remark goes for $G^\Omega_\omega$.
As mentioned in the introduction, the idea at the basis of our analysis is to somehow reduce the problem, 
locally, to a one-dimensional problem and to use
bounds that we have on the one-dimensional Green function in order to 
get exponential decay of the Green function for 
$\Gamma$. This is possible because of the fact that there is essentially only 
one path between any two vertices of the tree, and
because for any $\gamma>1$ one can find one-dimensional stretches of 
arbitrary length in $\Gamma$.

Let $x$ and $y$ be two distinct vertices of $\Gamma$. Then, by the resolvent formula,
\beq \label{resolvent1}
\begin{split}
G_\omega(z)&=G_{\calL(x,y);\omega}(z)-G_{\calL(x,y);\omega}(z) 
\left(H_\omega-H_{\calL(x,y);\omega}\right)
G_\omega(z) \\
&=G_{\calL(x,y);\omega}(z)-G_{\calL(x,y);\omega}(z) 
\left(\Delta-\Delta_{\calL(x,y)}\right)
G_\omega(z)
\end{split}
\eeq 
which holds in this form since $H_\omega$ and $H_{\calL;\omega}$ have the 
same diagonal part. Writing
\beq \label{resolvent2}
G_\omega(z)=G_{\calL(x,y)^{++};\omega}(z)-G_{\omega}(z) 
\left(\Delta-\Delta_{\calL(x,y)^{++}}\right)
G_{\calL(x,y)^{++};\omega}(z)
\eeq
(recall that $\calL(x,y)^{++}$ is the region in $\Gamma$ 
comprised of vertices of distance at most $2$ from 
 $\calL(x,y)$) 
and plugging this into \eqref{resolvent1}, we get:
\beq \label{double-resolvent}
\begin{split}
G_\omega(z)&=G_{\calL(x,y);\omega}(z)-G_{\calL(x,y);\omega}(z) 
T_{\calL(x,y)}G_{\calL(x,y)^{++};\omega}(z) \\
& \quad +G_{\calL(x,y);\omega}(z) T_{\calL(x,y)}G_{\omega}(z) 
T_{\calL(x,y)^{++}}
G_{\calL(x,y)^{++};\omega}(z),
\end{split}
\eeq
where we write 
\beq \no
T_\Omega=\Delta-\Delta_\Omega
\eeq
for any region $\Omega$ in $\Gamma$.

Now assume $x,y,w \in \scrV(\Gamma)$ are such that $y$ is on $\calL(x,w)$ and $w$ is outside of $\calL(x,y)^{++}$. Then,
\beq \no
\ip{\delta_{x}}{G_{\calL(x,y);\omega}(z) \delta_{w}}=0 
\eeq
and 
\beq \no 
\ip{\delta_{x}}{G_{\calL(x,y);\omega}(z) 
T_{\calL(x,y)}G_{\calL(x,y)^{++};\omega}(z) \delta_{w}}=0,
\eeq
so
\beq \label{double-resolvent-elements}
\begin{split}
\ip{\delta_{x}}{G_\omega(z) \delta_{w}}
&=\ip{\delta_{x}}{G_{\calL(x,y);\omega}(z) 
T_{\calL(x,y)}G_{\omega}(z) T_{\calL(x,y)^{++}}
G_{\calL(x,y)^{++};\omega}(z) \delta_{w}}\\
&=\sum_{(u,u') \in \Theta  (\calL(x,y) )} \sum_{(v,v') \in \Theta 
(\calL(x,y)^{++})}
\ip{\delta_{x}}{G_{\calL(x,y);\omega}(z) \delta_u} \\
& \quad \times \ip{\delta_{u'}}{G_{\omega}(z) \delta_v}
\ip{\delta_{v'}}{G_{\calL(x,y)^{++};\omega}(z)\delta_{w}}.
\end{split}
\eeq
It follows that
\beq \label{double-resolvent-s-mean}
\begin{split}
& \left \langle \left| \ip{\delta_{x}}{G_\omega(z) \delta_{w}} \right| ^s \right \rangle \\ 
&\leq \sum_{(u,u') \in \Theta  (\calL(x,y) )} \sum_{(v,v') \in \Theta 
(\calL(x,y)^{++})}
\Big \langle \left| \ip{\delta_{x}}{G_{\calL(x,y);\omega}(z) \delta_u} \right|^s  \\
& \times \left| \ip{\delta_{u'}}{G_{\omega}(z) \delta_v} \right|^s
 \left| \ip{\delta_{v'}}{G_{\calL(x,y)^{++};\omega}(z) \delta_{w}} \right|^s \Big \rangle
\end{split}
\eeq
for any $s \in (0,1)$ (we use $\langle \cdot \rangle$ to denote the mean over the disorder). As in \cite{aiz}, variants of 
equations 
\eqref{double-resolvent-elements} and 
\eqref{double-resolvent-s-mean} are the starting point of our derivation. 
We shall want to focus on the Green function restricted to some finite (large) balls 
$B(r)$ and not on the complete function. This is
no severe limitation as long as our estimates are uniform in $r$. We shall use the abbreviations
\beq \no
G^r_\omega \equiv G^{B(r)}_\omega,
\eeq
and 
\beq \no
T^r_\Omega=\Delta^{B(r)}-\Delta^{B(r)}_\Omega.
\eeq
Since the derivation of the first part of  
\eqref{double-resolvent-elements} uses only the resolvent formula, it is valid 
for $G^r$, so we have
\beq \label{restricted-double-resolvent-elements}
\begin{split}
&\ip{\delta_{x}}{G^r_\omega(z) \delta_{w}} 
=\ip{\delta_{x}}{G^r_{\calL(x,y);\omega}(z) 
T^r_{\calL(x,y)}G^r_{\omega}(z) T^r_{\calL(x,y)^{++}}
G^r_{\calL(x,y)^{++};\omega}(z) \delta_{w}} \\
&=\sum_{(u,u') \in \Theta  (\calL(x,y) )} \sum_{(v,v') \in \Theta 
(\calL(x,y)^{++})}
\ip{\delta_{x}}{G^r_{\calL(x,y);\omega}(z) \delta_u}
\ip{\delta_{u}}{T^r_{\calL(x,y)}\delta_{u'}} \\
& \quad \ip{\delta_{u'}}{G^r_{\omega}(z) \delta_v}
\ip{\delta_v}{T^r_{\calL(x,y)^{++}} \delta_{v'}}
\ip{\delta_{v'}}{G^r_{\calL(x,y)^{++};\omega}(z) \delta_{w}}.
\end{split}
\eeq
Here, $\ip{\delta_{u}}{T^r_{\calL(x,y)} \delta_{u'}} \leq 1$ since it may vanish if 
$u' \ {\rm or} \ u \notin B(r)$, and the same goes for $\ip{\delta_v}{T^r_{\calL(x,y)^{++}} \delta_{v'}}$. We see that 
\beq \label{restricted-double-resolvent-s-mean}
\begin{split}
& \left \langle \left| \ip{\delta_{x}}{G^r_\omega(z) \delta_{w}} \right|^s \right \rangle \\
&\leq \sum_{(u,u') \in \Theta^r (\calL(x,y) )} \sum_{(v,v') \in 
\Theta^r (\calL(x,y)^{++})}
\Big \langle \left| \ip{\delta_{x}}{G^r_{\calL(x,y);\omega}(z) \delta_u} \right| ^s \\
& \quad \times \left| \ip{\delta_{u'}}{G^r_{\omega}(z) \delta_v} \right|^s
 \left| \ip{\delta_{v'}}{G^r_{\calL(x,y)^{++};\omega}(z) \delta_{x}} \right|^s \Big \rangle \\
&\leq \sum_{(u,u') \in \Theta (\calL(x,y) )} \sum_{(v,v') \in \Theta 
(\calL(x,y)^{++})}
\Big \langle \left| \ip{\delta_{x}}{G^r_{\calL(x,y);\omega}(z) \delta_u} \right|^s \\
& \quad \times \left|\ip{\delta_{u'}}{G^r_{\omega}(z) \delta_v} \right|^s
 \left| \ip{\delta_{v'}}{G^r_{\calL(x,y)^{++};\omega}(z) \delta_{w}} \right|^s \Big \rangle,
\end{split}
\eeq
where $\Theta^r(\Omega)=\Theta(\Omega) \setminus \Theta(B(r))$. 

In order to get a useful bound from 
\eqref{restricted-double-resolvent-s-mean} and its variants, we need an a priori bound on the finite-volume restriction of 
the one-dimensional Green function. This is supplied for us by the following proposition due to Minami \cite{minami} 
(also see \cite{ckm, vdk}):

\begin{proposition}[Minami, Proposition A.1 of \cite{minami}] \label{Min}
Let $H_{\omega}=\Delta_{\bbZ}+\lambda V_\omega$ be a random Schr\"odinger 
operator on $\bbZ$ \rm{ (}with $\Delta_{\bbZ}$ the Laplacian on $\bbZ$, $\lambda>0${\rm)} and let $H_{\calL;\omega}$ be 
the restriction of $H_\omega$ to $\calL \subset \bbZ$ defined as above. Let
$G_{\calL;\omega}(x,y;z)=\ip{\delta_x}{ \left( H_{\calL;\omega}-z \right)^{-1}\delta_y}$. Assume that the random potential
consists of i.i.d.\ random variables with a common distribution $d\rho$ 
that satisfies 
\eqref{Minami-condition-1}-\eqref{Minami-condition-3}. 
Then for any $E \in \bbR$, there are $s_0 \in (0,1)$, $C>0$, $m>0$ and 
$\varepsilon>0$ such that
\beq \label{Minami-conclusion-1}
\left \langle |G_\calL(x,y;z)|^{s_0} \right \rangle \leq Ce^{-m|x-y|}
\eeq
for any finite segment $\calL \subset \bbZ$, $x \in \calL$, $y \in 
\mathscr{B}(\calL)$ and 
\beq \no
z \in \{z \mid \Ima z >0 , |z-E|<\varepsilon \}.
\eeq
\end{proposition}

Another bound we will need is a bound on the conditional expectation of the Green function restricted to any finite 
volume. For this we will use Lemma B.1 of \cite{aiz}. Note that the condition $R_1(\tau)$ of that paper is satisfied 
(with $\tau=1$) by 
any probability distribution $d\rho$ satisfying \eqref{Minami-condition-1}-\eqref{Minami-condition-3}. 
Thus, the hypotheses of Lemma B.1 of \cite{aiz} are satisfied in our situation, and we get 
\begin{proposition}[Aizenman et al., Lemma B.1 in \cite{aiz}] 
\label{Aiz}
There exists $\kappa<\infty$ such that for any finite 
subset $\Omega$ of $\scrV(\Gamma)$, any 
$x,y \in \Omega$, any $z \in \bbC$ and any $s \in (0,1)$,
\beq \label{Aizenman-conclusion}
\left \langle \left \vert \ip{\delta_x}{G_{\Omega;\omega}(z) \delta_y} \right \vert^s  \Big | 
\{V(u)\}_{u \in \Omega \setminus \{x,y\}} \right \rangle \leq \frac 
{1}{1-s} \frac {(4\kappa)^{s}}{\lambda^s},
\eeq 

where $\left \langle \left \vert \ip{\delta_x}{G_{\Omega;\omega}(z) \delta_y} \right \vert^s \Big | 
\{V(u)\}_{u \in \Omega \setminus \{x,y\}} \right \rangle$ is the conditional 
expectation of 
$\left \vert \ip{\delta_x}{G_{\Omega;\omega}(z) \delta_y}\right \vert^s$ 
conditioned on the values of the potential at all sites other than $x$ and $y$. 
\end{proposition}

Theorem \ref{localization} is implied, via standard arguments, by the following
\begin{theorem} \label{green}
Fix $x_0 \in \scrV(\Gamma)$. Then for any $E \in \bbR$, there 
are $s_0$, $A>0$ $q>0$ and $\varepsilon>0$ such 
that
\beq \label{conclusion}
\left \langle |G^r(x_0,v;z)|^{s_0} \right \rangle \leq Ae^{-q d(x_0,v)} 
\eeq
for any $r>|x_0|$, $v \in \scrV(\Gamma)$ and 
\beq \label{condition-on-z}
z \in \{z \mid \mathscr{I}z>0, |z-E|<\varepsilon \}.
\eeq
\end{theorem}

\begin{proof}[Proof of Theorem \ref{localization}]
Given Theorem \ref{green} and the fact that $\Gamma$ has finite dimensions 
in the sense of \eqref{dimensions}, this is a 
straightforward application of the Simon-Wolff criterion \cite[Theorem 2]{simon-wolff} (recall $d\rho$ is absolutely continuous 
with respect to Lebesgue measure). 
\end{proof}

\begin{proof}[Proof of Theorem \ref{green}]
Fix $r>|x_0|$ and $E \in \bbR$. Let $\kappa$ be the constant from Proposition \ref{Aiz} and let  
\beq
C(s)=\max \left( \frac {1}{1-s} \frac {(4\kappa)^{s}}{\lambda^s},1 \right).
\eeq

It follows by Proposition \ref{Min} that there exist $s_0 \in (0,1)$, $\varepsilon>0$, and $L>0$, such that if 
$x,y \in \scrV(\Gamma)$ are such that $\calL(x,y)$ has no junctions in it and $d(x,y) \geq L$, we will have
\beq \label{smallness}
\left \langle |G^r_{\calL(x,y)}(x,y;z)|^{s_0} \right \rangle<\frac{1}{4 C(s_0)^2}
\eeq
for all 
\beq \no
z \in \{z \mid \mathscr{I}z>0, |z-E|<\varepsilon \}.
\eeq
This is true because $G_{\calL(x,y)}(x,y;z)$ is just the same as the 
one-dimensional Green function restricted to a finite
segment, the boundary points of which are simply $x$ and $y$. The same 
goes for $G^r$ if $x$ and $y$ are both either in or
out of $B(r)$, (otherwise $|G^r_{\calL(x,y)}(x,y;z)|^{s_0}=0$ so \eqref{smallness} is still true). 

Fix
\beq \no
z \in \{z \mid \mathscr{I}z>0, |z-E|<\varepsilon \}.
\eeq
We want to iterate \eqref{restricted-double-resolvent-s-mean} with 
\eqref{smallness} in order to get the exponential decay we are trying to prove. 

Let $L_0=L+5$ and choose $n_0$ large enough so that $[\gamma^{n_0}]>8L_0$. Also, choose $R_0>|x_0|$ such that  
$R_0=\sum_{j=1}^{n_1}[\gamma^j]$ for some $n_1\geq n_0$ and $R_1=\sum_{j=1}^{n_1+1}[\gamma^j]$ 
and assume that $|v| \geq R_1$ (for the finite number of vertices $v \in B(R_1)$ we will bound the Green function 
by a constant). We may also assume that $r \geq |v|$, since, otherwise, $G^r(x_0,v;z)=0$. Let $u_{R_0}$ be the unique 
vertex on $\calL(x_0,v)$ with $|u_{R_0}|=R_0$. Let $\calL_0=\calL(x_0,u_{R_0})$. Then we have 
(see \eqref{restricted-double-resolvent-s-mean})
\beq \label{iteration0}
\begin{split}
\left \langle |\ip{\delta_{x_0}}{G^r_\omega(z) \delta_{v}}|^{s_0} \right \rangle &\leq 
\sum_{(u,u') \in \Theta (\calL_0 )} \sum_{(y,y') \in \Theta 
(\calL_0^{++})}
\Big \langle |\ip{\delta_{x_0}}{G^r_{\calL_0;\omega}(z) \delta_u}|^{s_0} \\
& \times \quad |\ip{\delta_{u'}}{G^r_{\omega}(z) \delta_y }|^{s_0}
 \left |\ip{\delta_{y'}}{G^r_{\calL_0^{++};\omega}(z) \delta_{v}} \right |^{s_0} \Big \rangle.
\end{split}
\eeq

Note that $G^r_{\calL_0^{++};\omega}(z)$ is a direct sum of operators, 
one corresponding to a finite tree containing $x_0$ and the others corresponding to various (infinite) forward trees 
emanating from the boundary points of $\calL_0^{++}$. 
Only one such tree contains $v$ so there is only one element $(y_0,y'_0) 
\in \Theta (\calL_0^{++})$ for which 
\beq \no
\left| \ip{\delta_{y'_0}}{G^r_{\calL_0^{++};\omega}(z) \delta_{v}} \right|^{s_0} \neq 0.
\eeq
Therefore
\beq \label{iteration0'}
\begin{split}
\left \langle \left| \ip{\delta_{x_0}}{G^r_\omega(z) \delta_{v}} \right|^{s_0} \right \rangle
&\leq \sum_{(u,u') \in \Theta (\calL_0 )} \Big \langle \left|\ip{\delta_{x_0}}
{G^r_{\calL_0;\omega}(z) \delta_u} \right|^{s_0} \\
& \quad \times \left|\ip{\delta_{u'}}{G^r_{\omega}(z) \delta_{y_0}} \right|^{s_0}
 \left|\ip{\delta_{y'_0}}{G^r_{\calL_0^{++};\omega}(z) \delta_{v}} \right|^{s_0} \Big \rangle.
\end{split}
\eeq

There are three terms on the RHS of the inequality above. The terms 
$\left|\ip{\delta_{x_0}}{G^r_{\calL_0;\omega}(z) \delta_u} \right|^{s_0}$ and
$\left|\ip{\delta_{y'_0}}{G^r_{\calL_0^{++};\omega}(z) \delta_{v}} \right|^{s_0}$
are independent random variables since the first depends only on the potential in $\calL_0$ and the second only 
on the potential outside of $\calL_0^{++}$. 
Moreover, neither of them depend on the potential at any of the $u'$ and at $y_0$. Therefore, we may evaluate 
the expectation by first evaluating the conditional expectation with respect to the potential at all other points. 
For this we may use Proposition \ref{Aiz} to get
\beq \label{iteration0''}
\begin{split}
& \left \langle \left|\ip{\delta_{x_0}}{G^r_\omega(z) \mid \delta_{v}} \right|^{s_0} \right \rangle \\ 
&\leq C(s_0) \sum_{u \in \mathscr{B}(\calL_0)} 
\left \langle \left|\ip{\delta_{x_0}}{G^r_{\calL_0;\omega}(z) \delta_u} \right|^{s_0} \right \rangle 
\left \langle \left|\ip{\delta_{y'_0}}{G^r_{\calL_0^{++};\omega}(z) \delta_{v}} \right|^{s_0} \right \rangle \\ 
&\leq C(s_0)^2 \left( \#B(R_0) \right)
\left \langle \left|\ip{\delta_{y'_0}}{G^r_{\calL_0^{++};\omega}(z) \delta_{v}} \right|^{s_0} \right \rangle.
\end{split}
\eeq

We proceed to estimate 
$\left \langle \left|\ip{\delta_{y'_0}}{G^r_{\calL_0^{++};\omega}(z) \delta_{v}} \right|^{s_0} \right \rangle$.
We start by dividing the line $\calL(y'_0,v)$ as follows: 
\begin{itemize}
\item Set $x_1=y'_0$.
Note that from our information about $v$ and the choice of $x_1$ (which 
reduces to the choice of $u_{R_0}$) $d(v,x_1)>7L_0$.
\item For any vertex $x \in \calL(x_1,v)$ let us denote by $\mathcal{J}(x)$ the 
distance from $x$ to the nearest junction on $\calL(x,v)$. Note that, by the choice of $R_0$ and 
since $L_0 \geq 5$, $\mathcal{J}(x_1)>L_0$.
Now, if $\mathcal{J}(x_1) \geq 3L_0$, let $v_1$ be the unique vertex at a 
distance $L_0$
from $x_1$ in $\calL(x_1,v)$. Otherwise, let $v_1$ be the unique vertex at 
a distance $5L_0$ from $x_1$ in $\calL(x_1,v)$.
\item Proceed by induction according to the following rule: Having defined 
$v_j$ for $j \geq 1$, let $x_{j+1}$ be the 
unique vertex at a distance $3$ from $v_j$ in $\calL(v_j,v)$. As long as 
$d(x_{j+1},v)>7L_0$, repeat the procedure above for 
choosing $v_{j+1}$, namely: If $\mathcal{J}(x_{j+1}) \geq 3L_0$, let $v_{j+1}$ be the unique 
vertex at a distance $L_0$ from $x_{j+1}$ in $\calL(x_{j+1},v)$. Otherwise, let $v_{j+1}$ be the unique vertex at a 
distance $5L_0$ from $x_{j+1}$ in $\calL(x_{j+1},v)$. If $d(x_{j+1},v) \leq 7L_0$, let $v_{j+1}=v$. 
\item We terminate the construction when $v_j=v$, of course. It is easy to 
see that this happens after a finite number of steps, since $d(x_j,x_{j+1})\leq 5L_0+3<7L_0$.
\end{itemize}

Thus, we get a set of vertex-pairs $\{(x_j,v_j)\}_{j=1}^l$ that satisfy:
\begin{enumerate}
\item For any $j$, $\mathcal{J}(x_j) \geq L_0$.
\item For any $j$, the distance between $v_j$ and the only junction (if there is one) on $\calL(x_j,v_j)$ is at 
least $L_0$.
\item For any $j$, $d(x_j,v_j) \geq L_0$.
\item $l \geq \frac{d(x_1,v)}{5L_0+3} \geq \frac{d(u_{R_0},v)}{6L_0}$.
\end{enumerate}

Let $\calL_j=\calL(x_j,y_j)$. We want to repeat the analysis leading to 
equations 
\eqref{restricted-double-resolvent-elements} and 
\eqref{restricted-double-resolvent-s-mean}, for 
$\left \langle \left|\ip{\delta_{y'_0}}{G^r_{\calL_0^{++};\omega}(z)\delta_{v}} \right|^{s_0} \right \rangle$.
We note that $G^r_{\calL_0^{++};\omega}(z)$ is the Green function of an 
operator for which the hopping terms have been 
removed both outside of $B(r)$ and for the boundary of $\calL_0^{++}$. 
Thus, setting 
\beq \no
\Omega_1=B(r) \cap \left( \Gamma \setminus \calL_0^{++} \right) 
\eeq 
and recalling that both $y'_0(=x_1)$ and $v$ are in $\left( \Gamma 
\setminus \calL_0^{++} \right)$, 
we get that 
\beq \no
\ip{\delta_{y'_0}}{G^r_{\calL_0^{++};\omega}(z)\delta_{v}}
=\ip{\delta_{x_1}}{G^{\Omega_1}_{\omega}(z)\delta_{v}}.  
\eeq
Note that $x_1$ is on the boundary of $\Omega_1$, so if we let 
\beq \no
T^{\Omega_1}_{\calL_1}=\Delta^{\Omega_1}-\Delta^{\Omega_1}_{\calL_1},
\eeq
we get that 
\beq \label{vanishing-boundary1'}
T^{\Omega_1}_{\calL_1} \delta_{x_1}=0.
\eeq

We have
\beq \label{iteration1}
\begin{split}
&\ip{\delta_{x_1}}{G^r_{\calL_0^{++};\omega}(z)\delta_{v}}
=\ip{\delta_{x_1}}{G^{\Omega_1}_{\omega}(z)\delta_{v}} \\ 
&=\sum_{(u,u') \in \Theta  (\calL_1 )} \sum_{(y,y') \in \Theta 
(\calL_1^{++})}
\ip{\delta_{x_1}}{G^{\Omega_1}_{\calL_1;\omega}(z)\delta_u}
\ip{\delta_{u}}{T^{\Omega_1}_{\calL_1}\delta_{u'}} 
\ip{\delta_{u'}}{G^{\Omega_1}_{\omega}(z)\delta_y} \\
& \quad \times \ip{\delta_{y}}{T^{\Omega_1}_{\calL_1^{++}}\delta_{y'}}
\ip{\delta_{y'}}{G^{\Omega_1}_{\calL_1^{++};\omega}(z)\delta_{v}}.
\end{split}
\eeq

Consider, first, $\calL_1^{++}$. As before, there is only one element 
$(y,y') \in \Theta (\calL_1^{++})$ for which 
$\ip{\delta_{y'}}{G^{\Omega_1}_{\calL_1^{++};\omega}(z) \delta_{v}} \neq 0$. From the construction,
it follows that this element is precisely $(x'_2,x_2)$ where $x'_2$ is the 
only point for which
$\ip{\delta_{x'_2}}{T^{\Omega_1}_{\calL_1^{++}}\delta_{x_2}} \neq 0.$ So
\beq \label{iteration1'}
\begin{split}
\ip{\delta_{x_1}}{G^r_{\calL_0^{++};\omega}(z)\delta_{v}} &=\sum_{(u,u') \in \Theta  (\calL_1 )}
\ip{\delta_{x_1}}{G^{\Omega_1}_{\calL_1;\omega}(z)\delta_u}
\ip{\delta_{u}}{T^{\Omega_1}_{\calL_1}\delta_{u'}} \\
& \quad \times \ip{\delta_{u'}}{G^{\Omega_1}_{\omega}(z)\delta_{x'_2}}
\ip{\delta_{x_2}}{G^{\Omega_1}_{\calL_1^{++};\omega}(z)\delta_{v}}.
\end{split}
\eeq

Now consider $\calL_1$. This is a linear path which has at most three points on its boundary. One is $x_1$, another is 
$v_1$. If $\calL_1$ has a junction in it (there is at most one in any case), then this is a third point on its boundary. 
There are no more possibilities. Because of \eqref{vanishing-boundary1'}, we have that the 
term corresponding to $x_1$ in the sum vanishes, so there are at most two terms 
in the sum above. Taking the $s_0$-moment for each of these terms and using Proposition \ref{Aiz} (by first averaging over 
the potential at $x'_2$ and $u'$, precisely as before) we get
\beq \label{iteration1''}
\begin{split}
\left \langle \left|\ip{\delta_{x_1}}{G^r_{\calL_0^{++};\omega}(z)\delta_{v}} \right|^{s_0} \right \rangle 
& \leq C(s_0) \sum_{u \in \mathscr{B} (\calL_1 ),\ u \neq x_1}
\left \langle \left|\ip{\delta_{x_1}}{G^{\Omega_1}_{\calL_1;\omega}(z)\delta_u} \right|^{s_0} \right \rangle \\
& \quad \times \left \langle \left|\ip{\delta_{x_2}}{G^{\Omega_1}_{\calL_1^{++};\omega}(z)\delta_{v}} \right|^{s_0} 
\right \rangle.
\end{split}
\eeq

We want to use \eqref{smallness} to estimate 
$\sum_{u \in \mathscr{B} (\calL_1 ),\ u \neq x_1} \left \langle \left|\ip{\delta_{x_1}}
{G^{\Omega_1}_{\calL_1;\omega}(z)\delta_u} \right|^{s_0} \right \rangle$.
Indeed, if $\calL_1$ contains no junctions then this sum has only one element, 
$\left \langle \left|\ip{\delta_{x_1}}
{G^{\Omega_1}_{\calL_1;\omega}(z)\delta_{v_1}} \right|^{s_0} \right \rangle$, and since it holds that 
\beq \label{restriction-equality}
\ip{\delta_{x_1}}{G^{\Omega_1}_{\calL_1;\omega}(z) \delta_u}
=\ip{\delta_{x_1}}{G^r_{\calL_1;\omega}(z) \delta_u},
\eeq
for any $u \in \scrV(\calL_1)$,
it immediately follows that 
\beq \label{smallness1}
\left \langle \left| \ip{\delta_{x_1}}{G^{\Omega_1}_{\calL_1;\omega}(z) \delta_{v_1}} \right|^{s_0} \right \rangle
< \frac{1}{4 C(s_0)^2}<\frac{1}{2C(s_0)}.
\eeq

If $\calL_1$ contains also a junction, $u_1$, then such a bound is not immediate from \eqref{smallness}. This case has 
two terms appearing in the sum: $\left \langle \left|\ip{\delta_{x_1}}
{G^{\Omega_1}_{\calL_1;\omega}(z)\delta_{u_1}} \right|^{s_0} \right \rangle$, and 
$\left \langle \left|\ip{\delta_{x_1}}
{G^{\Omega_1}_{\calL_1;\omega}(z)\delta_{v_1}} \right|^{s_0} \right \rangle$. Note that, since 
$\scrV(\calL_1) \subseteq \scrV(\Omega_1)$, we have that 
$\ip{\delta_{x_1}}{G^{\Omega_1}_{\calL_1;\omega}(z)\delta_{u_1}}=
\ip{\delta_{x_1}}{G^{\calL_1}_{\omega}(z)\delta_{u_1}}$, and  
$\ip{\delta_{x_1}}{G^{\Omega_1}_{\calL_1;\omega}(z)\delta_{v_1}}=
\ip{\delta_{x_1}}{G^{\calL_1}_{\omega}(z)\delta_{v_1}}$.

Let $u'_1$ be the unique backward neighbor of $u_1$ and let $\calL'_1=\calL(x_1,u'_1)$.
Consider, first 
$\left \langle \left|\ip{\delta_{x_1}}{G^{\calL_1}_{\omega}(z)\delta_{u_1}} \right|^{s_0} \right \rangle$. 
From \eqref{resolvent1} applied to $G^{\calL_1}_\omega(z)$, we get that 
\beq \label{resolvent1'}
\ip{\delta_{x_1}}{G^{\calL_1}_\omega(z)\delta_{u_1}}=-\ip{\delta_{x_1}}{G^{\calL_1}_{\calL'_1;\omega}(z)\delta_{u'_1}} 
\ip{\delta_{u_1}}{G_\omega(z)\delta_{u_1}}
\eeq
so, performing first the average over $V_\omega(u_1)$ (of which 
$\ip{\delta_{x_1}}{G^{\calL_1}_{\calL'_1;\omega}(z)\delta_{u'_1}}$ is independent), we get
\beq \label{smallness''}
\left \langle \left|\ip{\delta_{x_1}}
{G^{\Omega_1}_{\calL_1;\omega}(z)\delta_{u_1}} \right|^{s_0} \right \rangle =
\left \langle \left|\ip{\delta_{x_1}}
{G^{\calL_1}_{\omega}(z)\delta_{u_1}} \right|^{s_0} \right \rangle \leq \frac{C(s_0)}{4 C(s_0)^2}=\frac{1}{4 C(s_0)}
\eeq
(recall that $L_0=L+5$ so $d(x_1,u'_1)>L$).

As for $\left \langle \left|\ip{\delta_{x_1}}
{G^{\Omega_1}_{\calL_1;\omega}(z)\delta_{v_1}} \right|^{s_0} \right \rangle$, 
restricting $G^{\calL_1}_\omega(z)$ to $\calL'_1=\calL(x_1,u'_1)$ again, applying \eqref{double-resolvent}, taking the mean 
of the fractional moment and using Proposition \ref{Aiz}, we get that 
\beq \label{double-resolvent-s-mean2}
\begin{split}
\left \langle \left| \ip{\delta_{x_1}}{G^{\calL_1}_\omega(z) \delta_{v_1}} \right| ^s \right \rangle 
& \leq C(s_0) \left \langle \left| \ip{\delta_{x_1}}{G^{\calL_1}_{\calL'_1;\omega}(z) \delta_{u'_1}} \right|^s \right 
\rangle \\ 
& \quad \times \left \langle \left| \ip{\delta_{y'}}{G^{\calL_1}_{\calL^{'++}_1;\omega}(z) \delta_{v_1}} \right|^s \right \rangle
\end{split}
\eeq
where $y'$ is the only vertex on $\mathscr{B}(\calL^{'++}_1) \cap \calL_1$. Since $d(u_1,v_1)\geq L_0$ (see property 2 of the vertex 
pairs $\{x_j,v_j\}$) and $d(u_1,y')=2$, we have that $d(y',v_1) \geq L_0-2>L$. Furthermore, neither $\calL'_1$ nor 
$\calL(y',v_1)$ contains a junction, so the bound $\frac{1}{4 C(s_0)^2}$ applies to both 
Green functions on the RHS of \eqref{double-resolvent-s-mean2}. Thus
\beq \label{smallness'}
\left \langle \left| \ip{\delta_{x_1}}{G^{\calL_1}_\omega(z) \delta_{v_1}} \right| ^s \right \rangle 
\leq \frac{1}{16 C(s_0)^3}< \frac{1}{4 C(s_0)}.
\eeq
\eqref{smallness''} and \eqref{smallness'} give 
\beq \label{smallness2}
\sum_{u \in \mathscr{B} (\calL_1 ),\ u \neq x_1} \left \langle \left|\ip{\delta_{x_1}}
{G^{\Omega_1}_{\calL_1;\omega}(z)\delta_u} \right|^{s_0} \right \rangle \leq \frac{1}{2 C(s_0)}.
\eeq

Combining \eqref{smallness1} and \eqref{smallness2} with \eqref{iteration1''} 
we get 
\beq \label{base}
\left \langle \left|\ip{\delta_{x_1}}{G^r_{\calL_0^{++};\omega}(z) \delta_{v}} \right|^{s_0} \right \rangle 
\leq \frac{1}{2}
\left \langle \left|\ip{\delta_{x_2}}{G^{\Omega_1}_{\calL_1^{++};\omega}(z) \delta_{v}} \right|^{s_0} 
\right \rangle.
\eeq

At this point, we note that we can repeat the procedure outlined above for 
$\left \langle \left|\ip{\delta_{x_2}}{G^{\Omega_1}_{\calL_1^{++};\omega}(z) \delta_{v}} \right|^{s_0} 
\right \rangle$.
Writing 
\beq \no
\Omega_2=\Omega _1 \cap \left( \Gamma \setminus \calL_1^{++} \right)
\eeq 
we note that, as before, 
\beq \no
\ip{\delta_{x_2}}{G^{\Omega_1}_{\calL_1^{++};\omega}(z) \delta_{v}}
=\ip{\delta_{x_2}}{G^{\Omega_2}_{\omega}(z) \delta_{v}}
\eeq 
and also
\beq \no
T^{\Omega_2}_{\calL_2} \delta_{x_2}=0
\eeq
where
\beq \no
T^{\Omega_2}_{\calL_2}=\Delta^{\Omega_2}-\Delta^{\Omega_2}_{\calL_2}.
\eeq

Thus, we repeat the argument above with $\calL_2$ replacing $\calL_1$ and 
$x_3$ replacing $x_2$, to get the same estimate with 
$\left \langle \left|\ip{\delta_{x_2}}{G^{\Omega_1}_{\calL_1^{++};\omega}(z) \delta_{v}} \right|^{s_0} 
\right \rangle$
replaced by 
$\left \langle \left|\ip{\delta_{x_3}}{G^{\Omega_2}_{\calL_2^{++};\omega}(z) \delta_{v}} \right|^{s_0} 
\right \rangle.$

This can be repeated $l-1$ times, so that, estimating
\beq \no
\left \langle \left|\ip{\delta_{x_l}}{G^{\Omega_{l-1}}_{\calL_{l-1}^{++};\omega}(z) \delta_{v}} \right|^{s_0} \right \rangle
\leq C(s_0),
\eeq
we get (recall \eqref{iteration0''})
\beq \label{almost-conclusion}
\left \langle \left|\ip{\delta_{x_0}}{G^r_\omega(z) \delta_{v}} \right|^{s_0} \right \rangle
\leq C(s_0)^3 \left( \#B(R_0) \right) \left( \frac{1}{2} \right)^{l-1}. 
\eeq
$d(u_{R_0},v)=d(x_0,v)-d(u_{R_0},x_0)\geq d(x_0,v)-2R_0$ implies
\beq \label{almost-conclusion1}
\left \langle \left|\ip{\delta_{x_0}}{G^r_\omega(z) \delta_{v}} \right|^{s_0} \right \rangle
\leq A_1 e^{-q d(x_0,v)}, 
\eeq
for any $v$ with $|v| \geq R_1$, where $q = \frac{\ln 2}{6L_0}$ and 
$A_1=C(s_0)^3 \left( \#B(R_0) \right) 2^{(\frac{2R_0}{6L_0}+1)}$. Since there are only finitely many vertices in $B(R_1)$, 
it is obvious that one may choose a constant, $A>0$ so that
\beq \label{conclusion'}
\left \langle \left|\ip{\delta_{x_0}}{G^r_\omega(z) \delta_{v}} \right|^{s_0} \right \rangle
\leq A e^{-q d(x_0,v)}, 
\eeq
holds for any $v \in \scrV(\Gamma)$. Since none of our constants depended on $r$ or $z$, this estimate is uniform in 
$r>0$ and $z$ (in a proper neighborhood of $E$) so the conclusion follows.
\end{proof}



\begin{thebibliography}{99}

\bibitem{aiz}
M.~Aizenman, J.~H.~Schenker, R.~M.~Frierich and D.~Hundertmark, {\it Finite-volume fractional-moment criteria for Anderson 
localization}, Commun. Math. Phys. {\bf 224} (2001), 219--253.

\bibitem{aiz1} M.~Aizenman, R.~Sims and S.~Warzel, {\it Stability of the absolutely continuous spectrum of random 
Schr\"odinger operators on tree graphs}, Prob. Theor. Rel. Fields, to appear.

\bibitem{breuer-CMP} J.~Breuer, {\it Singular continuous spectrum for the Laplacian on certain sparse trees}, Commun. Math. Phys., to 
appear.

\bibitem{breuer-mol} J.~Breuer, {\it Singular continuous and dense point spectrum for sparse trees with finite dimensions}, 
Proceedings of ``Probability and Mathematical Physics'' a conference in honor of Stanislav Molchanov's 65th birthday, to appear. 

\bibitem{ckm} R.~Carmona, A.~Klein and F.~Martinelli, {\it Anderson localization for Bernoulli and other singular 
potentials}, Commun. Math. Phys. {\bf 108} (1987), 41--66.

\bibitem{carmona} R.~Carmona and J.~Lacroix, {\it Spectral Theory of Random Schr\"odinger Operators}, Birkh\"auser, Boston, 
1990.

\bibitem{froese-bethe} R.~Froese, D.~Hasler and W.~Spitzer, {\it Absolutely continuous spectrum for the Anderson model on a 
tree: a geometric proof of Klein's theorem}, Commun. Math. Phys., to appear.

\bibitem{klein-bethe} A.~Klein, {\it Extended states in the Anderson model on the Bethe lattice}, Adv. Math. {\bf 133} 
(1998), 163--184.

\bibitem{minami} 
N.~Minami, {\it Local fluctuation of the spectrum of a multidimensional Anderson tight binding model}, Commun. Math. Phys. 
{\bf 177} (1996), 709--725. 

\bibitem{simon-wolff} B.~Simon and T.~Wolff, {\it Singular continuous spectrum under rank one perturbations and 
localization for random Hamiltonians}, Commun. Pure Appl. Math. {\bf 39} (1986), no. 1, 75--90. 

\bibitem{vdk} H.~von~Dreifus and A.~Klein, {\it A new proof of localization in the Anderson tight binding model}, 
Commun. Math. Phys. {\bf 124} (1989), 285--299.

\end{thebibliography}
\end{document}